\documentclass[10pt]{article}
\font\smallit=cmti10
\font\smalltt=cmtt10

\usepackage{caption}
\usepackage{color}
\usepackage{amssymb,latexsym,amsmath,epsfig,amsthm} 
\usepackage{empheq}
\usepackage{here}
\usepackage{url}
\usepackage{ascmac}
\makeatletter
\usepackage{here}
\usepackage{comment}

\renewcommand\section{\@startsection {section}{1}{\z@}
{-30pt \@plus -1ex \@minus -.2ex}
{2.3ex \@plus.2ex}
{\normalfont\normalsize\bfseries\boldmath}}

\renewcommand\subsection{\@startsection{subsection}{2}{\z@}
{-3.25ex\@plus -1ex \@minus -.2ex}
{1.5ex \@plus .2ex}
{\normalfont\normalsize\bfseries\boldmath}}

\renewcommand{\@seccntformat}[1]{\csname the#1\endcsname. }

\makeatother
\newtheorem{theorem}{Theorem}
\newtheorem{lemma}{Lemma}
\newtheorem{conjecture}{Conjecture}

\newtheorem{corollary}{Corollary}
\theoremstyle{definition}
\newtheorem{definition}{Definition}

\newtheorem{pict}{Figure}[section]

\begin{document}

\begin{center}
\uppercase{\bf Wythoff's Game with a Pass}
\vskip 20pt
{\bf Ryohei Miyadera }\\
{\smallit Keimei Gakuin Junior and High School, Kobe City, Japan}\\
{\tt runnerskg@gmail.com}
\vskip 10pt
{\bf Hikaru Manabe}\\
{\smallit Tuskuba University, Tsukuba City, Japan}\\
{\tt urakihebanam@gmail.com}
\vskip 10pt
{\bf Masanori Fukui}\\
{\smallit Iwate Prefectural University }\\
{\tt fukui\_m@iwate-pu.ac.jp}

\end{center}
\vskip 20pt
\centerline{\smallit Received: , Revised: , Accepted: , Published: } 
\vskip 30pt


\centerline{\bf Abstract}
\noindent
This paper describes Wythoff's game with a pass, which is a variant of the classical Wythoff's game. The classical form is played with two piles of stones, from which two players take turns
to remove stones from one or both piles. When removing stones from both piles, an equal number must be removed from each. The player who removes the last stone or stones is the winner.
 In Wythoff's game with a pass, we modify the standard rules to allow for a one-time pass, i.e., a pass move that may be used at most once in a game but not from a terminal position. Once either player has used the pass, it is no longer available.
We denote the position of the game by $(x,y,p)$, where $x,y$ are numbers of stones in two piles and $p=1$ if a pass is available, and $p=0$ if not.
The authors proved that for $(x,y,1)$ with $x \geq 9$ or $y \geq 9$, $(x,y,1)$ is a P-position (the previous player's winning position) if and only if the Grundy number of $(x,y,0)$ is $1$.
Therefore, by using the result by U. Blass and A.S. Fraenkel, the Euclid distance between each previous player's winning position in Wythoff's game with a pass and a nearby previous player's winning position in Wythoff's game without a pass is within $\sqrt{20}$.
 \pagestyle{myheadings} 
 \markright{\smalltt \hfill} 
 \thispagestyle{empty} 
 \baselineskip=12.875pt 
 \vskip 30pt

\section{Introduction}\label{queenwithapass} 
Let  $\mathbb{Z}_{\geq0}$ and $\mathbb{N}$ be the sets of non-negative integers and natural numbers, respectively.
An interesting but challenging question in combinatorial game theory has been determining what happens when standard game rules are modified to allow a \textit{one-time pass}. This pass move may be used at most once in the game and not from a terminal position. Once either player has used a pass, it is no longer available. In the case of classical Nim, the introduction of the pass alters the mathematical structure of the game, considerably increasing its complexity. The effect of a pass on classical Nim remains an important open question that has defied traditional approaches. The late mathematician David Gale offered a monetary prize to the first person to develop a solution for three-pile classical Nim with a pass.

In \cite{nimpass} (p. 370), Friedman and Landsberg conjectured that ``solvable combinatorial games are structurally unstable to perturbations, while generic, complex games will be structurally stable.''One way to introduce such a perturbation is to allow a pass. 
One of the authors of the present article reported counterexamples to this conjecture in \cite{integers1} and \cite{integer2023}. The games used in \cite{integers1} and \cite{integer2023} are solvable because there are simple formulas for the Grundy numbers, and even when we introduce a pass move to the games, there are simple formulas for $\mathcal{P}$-positions.
Based on the research \cite{integers1} and \cite{integer2023}, the authors of the present article made the following conjecture.

\begin{conjecture}\label{conject}
Some games have specific mathematical structures that prevent the perturbation caused by the pass from spreading to other positions, and these games are solvable and have formulas for $\mathcal{P}$-positions even if a pass is introduced.
\end{conjecture}

Here, we present the research of Wythoff's game with a pass. Wythoff's game with a pass is not an exact example of Conjecture \ref{conject} because we do not have a formula for $\mathcal{P}$-positions. Still, it presents a very good example of specific mathematical structures that prevent the perturbation caused by the pass from spreading to other positions. 

For completeness, we briefly review some of the necessary concepts in combinatorial game theory by referring to $\cite{lesson}$ and \cite{combysiegel}.

\begin{definition}\label{NPpositions}
$ \mathrm{(i)}$	 A position is referred to as a $\mathcal{P}$-\textit{position} if it is the winning position for the previous player (the player who has just moved), as long as they play correctly at each stage. \\
$ \mathrm{(ii)}$ A position is referred to as an $\mathcal{N}$-\textit{position} if it is the winning position for the next player, as long as they play correctly at each stage.
\end{definition}
\begin{definition}\label{defofmexgrundy}
$ \mathrm{(i)}$	 For any position $\mathbf{p}$ of game $\mathbf{G}$, there is a set of positions that can be reached by precisely one move in $\mathbf{G}$, which we denote as \textit{move}$(\mathbf{p})$. \\	
$ \mathrm{(ii)}$	 The \textit{minimum excluded value} $(\textit{mex})$ of a set $S$ of non-negative integers is the least non-negative integer that is not in $S$. \\
$ \mathrm{(iii)}$	 Each position $\mathbf{p}$ of an impartial game $G$ has an associated \textit{Grundy number}, denoted by $\mathcal{G}(\mathbf{p})$.\\
The Grundy number is calculated recursively as: 
$\mathcal{G}(\mathbf{p}) = \textit{mex}\{\mathcal{G}(\mathbf{h}): \mathbf{h} \in \textit{move}(\mathbf{p})\}.$
\end{definition}

\begin{definition}\label{defofsum}
The \textit{disjunctive sum} of the two games, which is denoted as $\mathbf{G}+\mathbf{H}$, is a supergame in which a player may move in either $\mathbf{G}$ or $\mathbf{H}$ but not both.
\end{definition}

\begin{theorem}[\cite{lesson}]\label{theofgrundy}
$(i)$	 Let $\mathcal{G}$ be the Grundy number. Then $\mathbf{h}$ is a $\mathcal{P}$-position if and only if 
$\mathcal{G}(\mathbf{h})=0$.\\
$(ii)$	 The Grundy number of positions $\{\mathbf{g},\mathbf{h}\}$ in game $\mathbf{G}+\mathbf{H}$ is	$\mathcal{G}_{\mathbf{G}}(\mathbf{g})\oplus \mathcal{G}_{\mathbf{H}}(\mathbf{h})$, where $\oplus$ is nim-sum, or bitxor.
\end{theorem}
\section{Wythoff's Game}
\begin{figure}[H]
\begin{tabular}{cc}
\begin{minipage}[t]{0.5\textwidth}
\begin{center}
\includegraphics[height=3.cm]{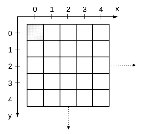}
\caption{    }
\label{chessboard}
\end{center}
\end{minipage}
\begin{minipage}[t]{0.5\textwidth}
\begin{center}
\includegraphics[height=2.5cm]{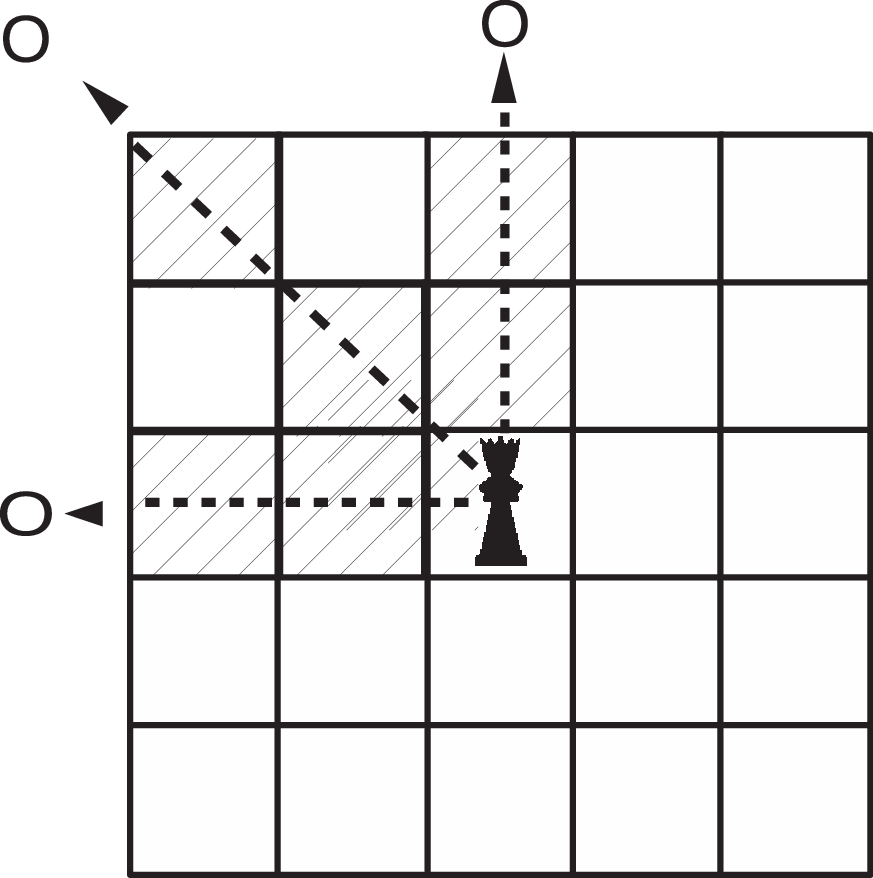}
\caption{   }
\label{moveofqueen}
\end{center}
\end{minipage}
\end{tabular}
\end{figure}

\begin{definition}\label{wythoff}
Wythoff's game is played with two piles of stones. Two players take turns removing stones from one or both piles; when removing stones from both piles, the number of stones removed from each must be equal. The player who removes the last stone or stones is the winner.
An equivalent description of the game is that a single chess Queen is placed somewhere on a large grid of squares, and each player can move the Queen towards the upper-left corner of the grid, either vertically, horizontally, or diagonally, any number of steps. 
The winner is the player who moves the Queen into the upper-left corner. See, Figure \ref{moveofqueen}.
\end{definition}

Figure~\ref{chessboard} shows the grid of squares, and we denote by $(x,y)$ the position of the Queen or the number of stones in the first and the second pile, where the horizontal coordinate and the vertical coordinate are denoted by 
$x$ and $y$.
Figure~\ref{moveofqueen} shows the moves the Queen can make in  Wythoff's game.
For the position $(x,y)$

\begin{theorem}[\cite{wythoffpaper}]\label{pofqueen}
The set of $\mathcal{P}$-positions of the game of Definition \ref{wythoff} is
	$\{ (\lfloor n \phi \rfloor, \lfloor n \phi \rfloor + n ):n \in \mathbb{Z}_{\ge 0} \}$
	$\cup \{ ( \lfloor n \phi \rfloor + n, \lfloor n \phi \rfloor ):n \in \mathbb{Z}_{\ge 0}  \}$,
 where $ \phi = \frac{1+\sqrt{5}}{2}$.
\end{theorem}

\begin{theorem}[\cite{fraenkelpaper}]\label{p1ofqueen}
Let $ \{(a_n,b_n):n \in \mathbb{Z}_{\ge 0}\}$ be the set of positions whose Grundy number are $1$.
Here, we assume that $a_n$ is increasing. Then,
	$\mid b_n - (\lfloor n \phi \rfloor + n) \mid \leq 4$ and 
	$ \lfloor n \phi \rfloor -1 \leq a_n \leq \lfloor n \phi \rfloor +2$.
\end{theorem}
This theorem is Corollary 5.14 of \cite{fraenkelpaper}.

\begin{corollary}\label{corowy}
For any position $(x,y)$ whose Grundy number is $1$, there exists a
position $(v,w)$ such that 
the Grundy number of $(v,w)$ is $0$ and the Euclid Distance of 
$(x,y)$ and $(v,w)$ is equal to or less than $\sqrt{20}$.
\end{corollary}
\begin{proof}
This is direct from Theorems \ref{pofqueen} and \ref{p1ofqueen}.   
\end{proof}

\section{Wythoff's Game with a Pass and the Sum of Wythoff's Game and a Pile of One Stone}
\begin{definition}\label{wythoffp}
We modify the standard rules of Wythoff's game so as to allow for a one-time pass. This pass move may only be used once in a game, not from a terminal position. Once either player has used the pass, it is no longer available.  
We denote by $\mathcal{G}_1$ the Grundy numbers of this game.
\end{definition}
\begin{definition}\label{wythoffs} 
By the sum of two games in Definition \ref{defofsum}, we define the sum of the classical Wythoff's game without a pass and the game of a pile of one stone. 
We denote by $\mathcal{G}_2$ the Grundy numbers of the game of this game.
\end{definition}

We denote the position of the Queen with three coordinates
$\{x,y,p\}$ for the game of Definition \ref{wythoffp} and the game of Definition \ref{wythoffs}.
The coordinates $x,y$ define the number of stones in the first and the second pile, or if we use Queen in the game, the position of the Queen on the chessboard. As for the game of Definition \ref{wythoffp}, 
the additional parameter $p$ denotes whether the pass is still available ($p=1$) or has already been used ($p=0$). 
For the game of Definition \ref{wythoffs}, the parameter $p=1$ if there is a stone in the third pile, and $p=0$ if not. Note that when $p=0$, the games of  Definitions \ref{wythoffp} and \ref{wythoffs} are the classical Wythoff's game.

\begin{definition}\label{movewythoff}
 For any $x,y \in \mathbb{Z}_{\ge 0}$ and $p=0,1$, let
\begin{equation}
M_1(x,y,p)= \{(u,y,p):u<x \text{ and } u \in \mathbb{Z}_{\ge 0}\},  \label{m1set}
\end{equation}

\begin{equation}
M_2(x,y,p)=\{(x,v,p):v<y \text{ and } v \in \mathbb{Z}_{\ge 0}\}, \label{m2set}  
\end{equation}

\begin{equation}
M_3(x,y,p)=\{(x-t,y-t,p): 1 \leq t \leq \min(x,y) \text{ and } t \in \mathbb{Z}_{\ge 0}\},\label{m3set}
\end{equation}

\begin{equation}
M_4(x,y,p)=
\begin{cases}
\{(x,y,0)\} & (\mbox{ if } x+y >0 \text{ and }p=1),\\
\emptyset & (\mbox{ if } x+y = 0 \text{ or }p=0 ),
\end{cases}\label{m4set}
\end{equation} 
and
\begin{equation}
M^{\prime}_4(x,y,p)=
\begin{cases}
\{(x,y,0)\} & (\mbox{ if } p=1),\\
\emptyset & (\mbox{ if }p=0 ).
\end{cases}\label{m4pset}
\end{equation} 
\end{definition}
Sets $M_1(x,y,p)$, $M_2(x,y,p)$, and $M_3(x,y,p)$ are the sets of horizontal, vertical, and diagonal moves, respectively.
Set $M_4(x,y,p)$ is the set of the pass move of Wythoff's game with a pass in Definition \ref{wythoffp}, and 
Set $M^{\prime}_4(x,y,p)$  is the set of moves in the third pile of the game of Definition \ref{wythoffs}.
Note that 
$M_4(x,y,p)$ is empty if and only if $x+y = 0$ or $p=0$, and $M^{\prime}_4(x,y,p)$ is empty if and only if $p=0$.

Next, we define  $move_1$ and $move_2$, which are moves of the games of Definitions \ref{wythoffp} and \ref{wythoffs}, respectively.
\begin{definition}\label{movewythoff2}
 For any $x,y \in \mathbb{Z}_{\ge 0}$ and $p=0,1$, let
\begin{equation}
move_1(x,y,p)=M_1(x,y,p) \cup M_2(x,y,p) \cup M_3(x,y,p) \cup M_4(x,y,p)    
\end{equation}
and 
\begin{equation}
move_2(x,y,p)=M_1(x,y,p) \cup M_2(x,y,p) \cup M_3(x,y,p) \cup M^{\prime}_4(x,y,p).   
\end{equation}
\end{definition}
\section{The positions $(x,y,p)$ such that $x,y \leq 8$}
This section aims to prove Lemma \ref{lemmaforabc}, but the readers of this article can skip the proof of Lemma \ref{lemmaforabc} and 
read the programs in Section \ref{computercalculation}, where three programs are presented. Since the conclusion of Lemma \ref{lemmaforabc} is about sets of a small number of positions, a reliable computer program can replace the proof of the lemma.

We now define three sets, $A$, $B$, and $C$, and study them. These sets have mathematical structures mentioned in Conjecture \ref{conject} that prevent the perturbation caused by the pass from spreading to other positions.

\begin{figure}[H]
\begin{tabular}{ccc}
\begin{minipage}[t]{0.33\textwidth}
\begin{center}
\includegraphics[height=3.4cm]{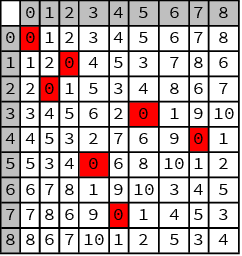}
\caption{Set $A$}
\label{nopass}
\end{center}
\end{minipage}
\begin{minipage}[t]{0.33\textwidth}
\begin{center}
\includegraphics[height=3.4cm]{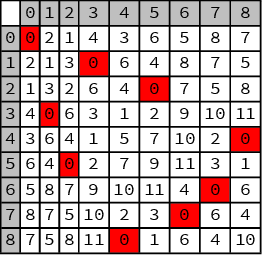}
\caption{Set $B$ }
\label{queenpass15}
\end{center}
\end{minipage}
\begin{minipage}[t]{0.33\textwidth}
\begin{center}
\includegraphics[height=3.4cm]{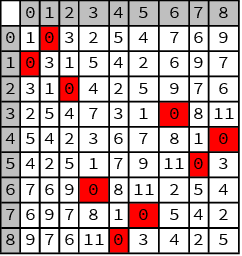}
\caption{Set $C$}
\label{queenplusone}   
\end{center}
\end{minipage}
\end{tabular}
\end{figure}

\begin{definition}\label{defnofab} 
	Let 
\begin{equation}
A=\{(0, 0, 0), (1, 2, 0), (2, 1, 0), (3, 5, 0), (4, 7, 0), (5, 3, 0), (7,
   4, 0)\}, \nonumber
\end{equation}
\begin{equation}
B= \{(0,0,1),(1,3,1),(3,1,1),(2,5,1),(5,2,1),(4,8,1),(8,4,1),(6,7,1),(7,6,1) \}, \nonumber
\end{equation}  
and 
\begin{equation}
C = \{(0,1,1),(1,0,1),(2,2,1),(3,6,1),(6,3,1),(4,8,1),(8,4,1),(5,7,1),(7,5,1) \}.\nonumber
\end{equation}  
\end{definition}
In Figures \ref{nopass}, \ref{queenpass15}, and \ref{queenplusone}, we have sets 
$A,B,$ and $C$. They are printed in red.
  
\begin{figure}[H]
\begin{tabular}{cc}
\begin{minipage}[t]{0.5\textwidth}
\begin{center}
\includegraphics[height=3.6cm]{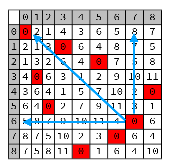}
\caption{From P not P}
\label{queenarrow}
\end{center}
\end{minipage}
\begin{minipage}[t]{0.5\textwidth}
\begin{center}
\includegraphics[height=3.4cm]{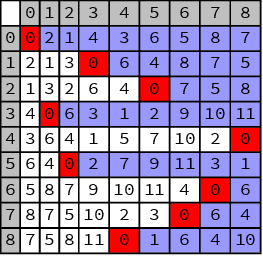}
\caption{horizontal }
\label{hmovep}
\end{center}
\end{minipage}
\end{tabular}
\end{figure}

\begin{figure}[H]
\begin{tabular}{cc}
\begin{minipage}[t]{0.5\textwidth}
\begin{center}
\includegraphics[height=3.4cm]{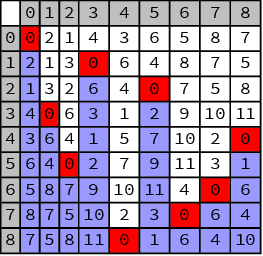}
\caption{vertical}
\label{vmovep}
\end{center}
\end{minipage}
\begin{minipage}[t]{0.5\textwidth}
\begin{center}
\includegraphics[height=3.4cm]{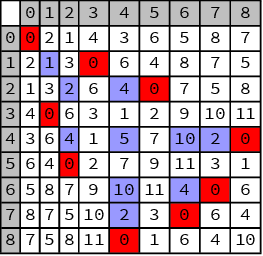}
\caption{daiagonal}
\label{dmovep}
\end{center}
\end{minipage}
\end{tabular}
\end{figure}

\begin{lemma}\label{lemmaforabc}
$(i)$
The set $A$ in Definition \ref{defnofab}  is the set of $\mathcal{P}$-positions $(x,y,0)$ of the game of  Definitions \ref{wythoffp} such that $x,y \leq 8$, where $x,y$ are horizontal and vertical coordinates when the pass is not available.\\
$(ii)$
The set $B$ in Definition \ref{defnofab} is the set of $\mathcal{P}$-positions  $(x,y,1)$ of the game of  Definitions \ref{wythoffp} such that $x,y \leq 8$ when the pass is available.\\
$(iii)$
The set $C$ in Definition \ref{defnofab}  is the set of $\mathcal{P}$-positions  $(x,y,1)$ of the game of  Definitions \ref{wythoffs}  such that $x,y \leq 8$ when the third coordinate is $1$.
\end{lemma}
\begin{proof}
By using Theorem \ref{pofqueen} for $x,y \leq 8$ we prove $ \mathrm{(i)}$.\\
$ \mathrm{(ii)}$	
Let $U=\{(x,y,1):x,y \leq 8\}$. Since we prove that the set $\mathcal{P}$-positions of the game of  Definitions \ref{wythoffp} in $U$ is $B$, we need to prove that
\begin{equation}
move_1(x,y,1) \cap (A \cup B) = \emptyset \label{moveemplty}
\end{equation}
for any $(x,y,1) \in B$
and 
\begin{equation}
move_1(x,y,0) \cap (A \cup B) \ne \emptyset\label{move1b} 
\end{equation}
for any $(x,y,1) \in U-B$.

First, we prove Relation (\ref{moveemplty}). Suppose that we start with the position $(7,6,1)$.
Then, the horizontal, the vertical and the diagonal moves are described in
Figure \ref{queenarrow}, and it is easy to see that 
$M_1(7,6,1) \cap B = \emptyset$, $M_2(7,6,1) \cap B = \emptyset$, and $M_3(7,6,1) \cap B = \emptyset$. $M_4(7,6,1)=(7,6,0) \notin A$.
Therefore, $move_1(7,6,1) \cap (A \cap B) = \emptyset$.

For  any $(x,y,1) \in B$, by Figure \ref{queenpass15}, it is easy to show that    
$M_1(x,y,1) \cap B = \emptyset$, $M_2(x,y,1) \cap B = \emptyset$, and 
$M_3(x,y,1) \cap B = \emptyset$. 
By comparing Figure \ref{nopass} and Figure \ref{queenpass15}, we have $M_4(x,y,1) \cap A = \emptyset$. 
Therefore, we have Relation \ref{moveemplty}.\\
Let  $(x,y,1) \in U-B$. We prove Relation \ref{move1b}.
For $(8,y,1)$ with $y=0,1,2,3$ and $y=5,6,7,8$, it is clear that $M_1(x,y,1) \cap B \ne \emptyset$.
In this way, for all the blue positions $(x,y,1)$ in Figure \ref{hmovep}, $M_1(x,y,0) \cap B \ne \emptyset$.
Similarly, we prove that for all the blue positions $(x,y,1)$ in Figure \ref{vmovep}, $M_2(x,y,0) \cap B \ne \emptyset$ and for all the blue positions $(x,y,1)$ in Figure \ref{dmovep}, $M_3(x,y,1) \cap B \ne \emptyset$.
The position in $U$ that does not belong to the blue positions in Figures \ref{hmovep}, \ref{vmovep}, and \ref{dmovep} is $(1,2,1)$ or $(2,1,1)$. Since $(1,2,0)$ and $(2,1,0)$ belong to the set $A$ in Figure \ref{nopass}.
$M_4(1,2,1)\cap A \ne \emptyset$ and $M_4(2,1,1)\cap A \ne \emptyset$. Therefore, we have Relation \ref{move1b}.\\
$ \mathrm{(iii)}$	
By the method that is very similar to the one used in $ \mathrm{(ii)}$, we can prove $ \mathrm{(iii)}$. Therefore, the details are omitted. 
\end{proof}

\section{The Set of $\mathcal{P}$-positions of Wythoff's Game with a Pass}

Let $P_0=\{(x,y,0):\mathcal{G}_1(x,y,0)=0\}$, $P_1=\{(x,y,1):\mathcal{G}_1(x,y,1)=0\}$, and $P_2=\{(x,y,0):\mathcal{G}_2(x,y,1)=0\}$.
\begin{lemma}\label{lemmamove}
Let $x,y \in \mathbb{Z}_{\geq0}$ such that $x \geq 9$ or $y \geq 9$. Then, we obtain the following.\\
$(i)$ If  $ x \geq 9$ and $y \leq 8$,  $M_1(x,y,1) \cap B \ne \emptyset$, $M_1(x,y,1) \cap C \ne \emptyset$, and 
$M_2(x,y,1) \cap C = M_2(x,y,1) \cap B = \emptyset$.\\
$(ii)$ If $x \leq 8$ and $y \geq 9$, $M_2(x,y,1) \cap B \ne \emptyset$, $M_2(x,y,1) \cap C \ne \emptyset$, and $M_1(x,y,1) \cap C = M_1(x,y,1) \cap B = \emptyset$.\\
$(iii)$ If $x \leq y+4$ or $y \leq x+4$, 
$M_3(x,y,1) \cap B  \ne \emptyset$ and 
$M_3(x,y,1) \cap C \ne \emptyset$.\\
$(iv)$ If $x \geq y+5$ or $y \geq x+5$, $M_3(x,y,1) \cap B = M_3(x,y,1) \cap C = \emptyset$.
\end{lemma}
\begin{proof}
$ \mathrm{(i)}$
By Definition \ref{defnofab}, there exist $v, v^{\prime} \in \mathbb{Z}_{\geq0}$ such that $1 \leq u,u^{\prime} \leq 8$, $(u,y,1) \in B$ and $(u^{\prime},y,1) \in C$. Then, 
$(u,y,1) \in M_1(x,y,1) \cap B$ and $(u^{\prime},y,1) \in M_1(x,y,1) \cap C$, and $M_1(u,v,1) \cap B \ne \emptyset$ and $M_1(u,v,1) \cap C \ne \emptyset$.
 Since $x \geq 9$, $M_2(x,y,1) \subset \{(x,v,1):v \in \mathbb{Z}_{\geq0} \}$ $ \subset (B \cup C)^c$, where $(B \cup C)^c$ is the complement of the set $B \cup C)$. Hence, $M_2(x,y,1) \cap C = M_2(x,y,1) \cap B = \emptyset$.\\
 $ \mathrm{(ii)}$ is direct from $ \mathrm{(i)}$, since this game is symmetrical with respect to the first and the second coordinates.\\
$ \mathrm{(iii)}$ We have two cases, 
\noindent {\tt Case 1}: Suppose that $x \leq y+4$. By Definition \ref{defnofab}, for any $a \in \mathbb{Z}_{\geq0}$ such that $a \leq 4$, there exist $u,u^{\prime},v,v^{\prime} \in \mathbb{Z}_{\geq0}$ 
such that $u=v+a$, $u^{\prime}=v^{\prime}+a$, $(u,v,1) \in B$, and $(u^{\prime},v^{\prime},1) \in C$.
Then, $(u,v,1) \in B \cup M_3(x,y,1)$ and $(u^{\prime},v^{\prime},1) \in C \cup M_3(x,y,1)$

Therefore, $M_3(x,y,1) \cap B  \ne \emptyset$ and $M_3(x,y,1) \cap C \ne \emptyset$.\\
\noindent {\tt Case 2}: Suppose that $y \leq x+4$.  Since this game is symmetrical with respect to the first and the second coordinates,
$M_3(x,y,1) \cap B  \ne \emptyset$ and $M_3(x,y,1) \cap C \ne \emptyset$.\\
$ \mathrm{(iv)}$ We have two cases.\\
\noindent {\tt Case 1}: Suppose that $x \geq y+5$. There is no $u,v \in \mathbb{Z}_{\geq0}$ such that $u \geq v+5$ and 
$(u,v) \in B \cup C$. Hence, $M_3(x,y,1) \cap B = M_3(x,y,1) \cap C = \emptyset$.\\
\noindent {\tt Case 2}: Suppose that $y \geq x+5$. Since this game is symmetrical with respect to the first and the second coordinates,
 $M_3(x,y,1) \cap B = M_3(x,y,1) \cap C = \emptyset$.
\end{proof}

\begin{lemma}\label{lemmamove2}
For $x,y$ such that $x \geq 9$ or $y \geq 9$, we have the following.\\ 
$(i)$ $M_1(x,y,1) \cap B \ne \emptyset$ if and only if $M_1(x,y,1) \cap C \ne \emptyset$.\\
$(ii)$ $M_2(x,y,1) \cap B \ne \emptyset$ if and only if $M_2(x,y,1) \cap C \ne \emptyset$.\\
$(iii)$ $M_3(x,y,1) \cap B \ne \emptyset$ if and only if $M_3(x,y,1) \cap C \ne \emptyset$.
\end{lemma}
\begin{proof}
By Lemma \ref{lemmamove}, we have $\mathrm{(i)}$, $\mathrm{(i)}$, and $\mathrm{(iii)}$. 
\end{proof}

\begin{theorem}\label{twogamesth}
For $x,y \in \mathbb{Z}_{\ge 0}$ such that $x \geq 9$ or $y \geq 9$,
\begin{equation}
\mathcal{G}_1(x,y,1)=0 \text{ if and only if } \mathcal{G}_2(x,y,1)=0.\label{g1g2eq}
\end{equation}
\end{theorem}
\begin{proof}
Let $V_8=\{(x,y,1):x,y \leq 8\}$.
It is sufficient to prove that
\begin{equation}
P_1 - V_8 =  P_2 - V_8.
\end{equation}
Let $U_k=\{(x,y,1):x+y \leq k\}$, and we prove that 
\begin{equation}
(U_{n}-V_8) \cap P_1 = 
(U_{n}-V_8) \cap P_2  
\end{equation}
for any natural number $n$ by mathematical induction.
Since $(U_{17}-V_8) \subset \{(u,v):u \geq 9 \text{ and } v \leq 8\}$
$\cup \{(u,v):u \leq 8 \text{ and } v \geq 9\}$, by  $(i)$ and $(ii)$ of Lemma \ref{lemmamove}, any point in $U_{17}-V_8$ is a $\mathcal{N}$-position. Hence 
\begin{equation}
(U_{17}-V_8) \cap P_1 = 
(U_{17}-V_8) \cap P_2 = \emptyset.\nonumber
\end{equation}
Suppose that for some natural number $k$
\begin{equation}
(U_{k}-V_8) \cap P_1 = 
(U_{k}-V_8) \cap P_2.   \label{ukv8requa} 
\end{equation}
Let  $x,y \in \mathbb{Z}_{\ge 0}$ such that $(x,y,1) \in U_{k+1}-V_8$.
Then, for $i=1,2,3$, by Definition \ref{movewythoff}
\begin{equation}
M_i(x,y,1) \subset U_k,
\end{equation}
and hence we have 
\begin{align}
M_i(x,y,1)\cap P_1& = M_i(x,y,1) \cap ((U_k-V_8) \cup V_8) \cap P_1 \nonumber \\
& = (M_i(x,y,1) \cap (U_k-V_8)\cap P_1) \cup (M_i(x,y,1) \cap V_8\cap P_1)  \nonumber \\
& = (M_i(x,y,1) \cap (U_k-V_8)\cap P_1) \cup (M_i(x,y,1) \cap B) \label{m1andb} 
\end{align}
and
\begin{align}
M_i(x,y,1)\cap P_2& = M_i(x,y,1) \cap ((U_k-V_8) \cup V_8) \cap P_2 \nonumber \\
& = (M_i(x,y,1) \cap (U_k-V_8)\cap P_2) \cup (M_i(x,y,1) \cap V_8\cap P_2)  \nonumber \\
& = (M_i(x,y,1) \cap (U_k-V_8)\cap P_2) \cup (M_i(x,y,1) \cap C). \label{m1andb2}
\end{align}      
By Lemma \ref{lemmamove2}, Equations  (\ref{ukv8requa}), (\ref{m1andb}), and (\ref{m1andb2}), we have 
\begin{equation}
M_i(x,y,1)\cap P_1 \ne \emptyset \text{ if and only if } M_i(x,y,1)\cap P_2 \ne \emptyset \label{mirelation}
\end{equation}
for $i=1,2,3$.

Since $M_4(x,y,1)=M_4^{\prime}(x,y,1)=(x,y,0)$ for $x,y$ such that $x \geq 9$ or $y \geq 9$, we have 
\begin{equation}
M_4(x,y,1) \cap P_0 \ne \emptyset \text{ if and only if } M^{\prime}_4(x,y,1) \cap P_0 \ne \emptyset.\label{m4relation}
\end{equation}

By (\ref{mirelation}) and (\ref{m4relation}), for any $x,y$ such that $(x,y,1) \in U_{k+1}-V_8$,
\begin{equation}
move_1(x,y,1)\cap (P_1\cup P_0) \ne \emptyset \text{ if and only if } move_2(x,y,1)\cap (P_2 \cup P_0) \ne \emptyset,\nonumber 
\end{equation}
and hence,
\begin{equation}
(U_{k+1}-V_8) \cap P_1 = 
(U_{k+1}-V_8) \cap P_2.   \nonumber
\end{equation}
Therefore, by mathematical induction, we have 
\begin{equation}
(U_{n}-V_8) \cap P_1 = 
(U_{n}-V_8) \cap P_2  \nonumber 
\end{equation}
 for any natural number $n$ 
\end{proof}

By Theorem \ref{twogamesth}, a position $(x,y,1)$ is a $\mathcal{P}$-position of the game in Definition \ref{wythoffp} if and only if it is a $\mathcal{P}$-position of the game in Definition \ref{wythoffs} when $x, y \geq 9$. 
The mathematical structure of Set $B$ prevents the perturbation caused by the pass from spreading to other positions.

\begin{lemma}\label{ififth}
For any $x,y \in \mathbb{Z}_{\ge 0}$ such that $x \geq 9$ or $y \geq 9$, 
\begin{equation}
\mathcal{G}_1(x,y,1)= 1 \text{ if and only if }  \mathcal{G}_1(x,y,0)= 0. \label{ifandonlyif}
\end{equation}
\end{lemma}
\begin{proof}
By Theorem \ref{theofgrundy} and Definition \ref{wythoffs},
$\mathcal{G}_1(x,y,1)=\mathcal{G}_1(x,y,0) \oplus 1$. Hence, we have
Relation (\ref{ifandonlyif}).
\end{proof}

\begin{theorem}\label{theofg1g0}
For $x,y \in \mathbb{Z}_{\ge 0}$ such that $x \geq 9$ or $y \geq 9$,
 $\mathcal{G}_1(x,y,1)=0$ if and only if  $\mathcal{G}_1(x,y,0)=1$.
\end{theorem}
\begin{proof}
This is direct from Theorem \ref{twogamesth} and Lemma \ref{ififth}.
\end{proof}

\begin{theorem}\label{g0g1the}
For any position $(x,y)$ such that 
 $\mathcal{G}_1(x,y,1)=0$, there exist  $(v,w)$ such that  $\mathcal{G}_1(x,y,0)=0$
 and the Euclid Distance of $(x,y)$ and $(v,w)$ is within $\sqrt{20}$.
\end{theorem}
\begin{proof}
This is direct from Corollary \ref{corowy} and Theorem \ref{theofg1g0}.
\end{proof}

By Theorem \ref{g0g1the}, the Graph of the set of $\mathcal{P}$-positions in the classical Wythoff's game and the graph of $\mathcal{P}$-positions in Wythoff's game with a pass when the pass is still available look very similar.  See Figures \ref{wythoffclassic} and \ref{wythoffpass1}.

\begin{center}
	\includegraphics[height=3.5cm]{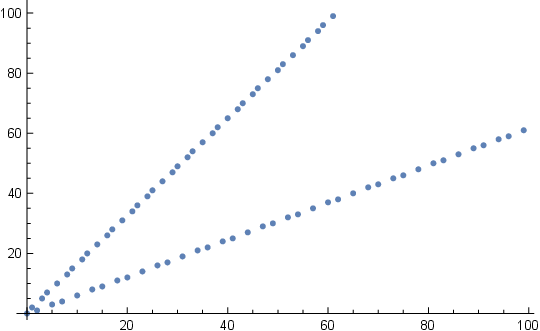}
\begin{pict}\label{wythoffclassic}
	The Graph of the set of $\mathcal{P}$-positions in the classical Wythoff's game.
	\end{pict}
	
	\includegraphics[height=3.5cm]{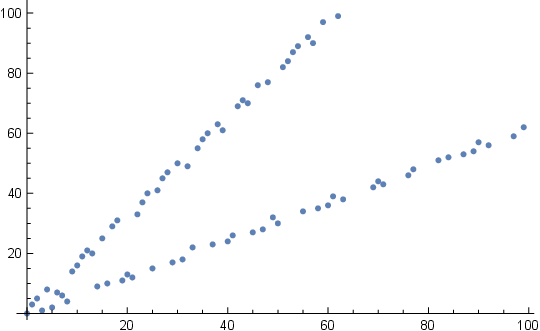}
	\begin{pict}\label{wythoffpass1}
The graph of $\mathcal{P}$-positions in Wythoff's game with a pass when the pass is still available.
	\end{pict}
\end{center}

\section{Computer Programs that Calculate the Set $A$,$B$, and$C$}\label{computercalculation}
Here, we have a Mathematica program, a CGSuite program, and a Python program to prove Lemma \ref{lemmaforabc}. 
These three programs were made by the three authors of the present article independently using three different computer languages.
Since the conclusion of Lemma \ref{lemmaforabc} is about sets of a small number of positions, a reliable computer program can replace the proof of the lemma.
\subsection{Calculation by Mathematica}
In this subsection, we give Tables of Grundy numbers using Mathematica.
First, we present a Mathematica program for calculating the Grundy numbers for Wythoff's game.
\begin{figure}[H]
\begin{tabular}{ccc}
\begin{minipage}[t]{0.33\textwidth}
\begin{center}
\includegraphics[height=3.4cm]{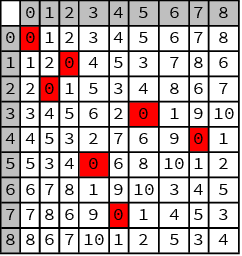}
\caption{Set $A$}
\label{nopasstable}
\end{center}
\end{minipage}
\begin{minipage}[t]{0.33\textwidth}
\begin{center}
\includegraphics[height=3.4cm]{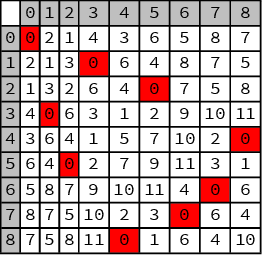}
\caption{Set $B$ }
\label{queenpasstable}
\end{center}
\end{minipage}
\begin{minipage}[t]{0.33\textwidth}
\begin{center}
\includegraphics[height=3.4cm]{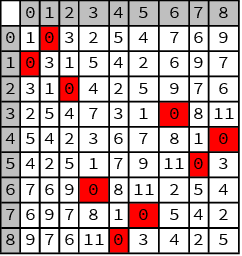}
\caption{Set $C$}
\label{queenplusonetable}   
\end{center}
\end{minipage}
\end{tabular}
\end{figure}

\begin{small}
\begin{verbatim}
ss = 15; allcases = 
Flatten[Table[{a, b, c}, {a, 0, ss}, {b, 0, ss}, {c, 0, 1}], 2];
move[z_] := Block[{p}, p = z; If[p[[1]] + p[[2]] > 0,
  Union[Table[{t1, p[[2]], p[[3]]}, {t1, 0, p[[1]] - 1}], 
   Table[{p[[1]], t2, p[[3]]}, {t2, 0, p[[2]] - 1}], 
   Table[{p[[1]], p[[2]], t3}, {t3, 0, p[[3]] - 1}],
   Table[{p[[1]] - tu, p[[2]] - tu, p[[3]]}, {tu, 1, Min[p[[1]], p[[2]]]}]], 
  Union[Table[{t1, p[[2]], p[[3]]}, {t1, 0, p[[1]] - 1}], 
   Table[{p[[1]], t2, p[[3]]}, {t2, 0, p[[2]] - 1}],
   Table[{p[[1]] - tu, p[[2]] - tu, p[[3]]}, {tu, 1, Min[p[[1]], p[[2]]]}]]]];
Mex[L_] := Min[Complement[Range[0, Length[L]], L]];
Gr[pos_] := Gr[pos] = Mex[Map[Gr, move[pos]]];
pposition = Select[allcases, Gr[#] == 0 &];
gg = 0.75; b2 = Table[{1, t}, {t, 2, ss + 2}]; 
b3 =  Table[{t, 1}, {t, 2, ss + 2}];
b1 = Table[{t, 9}, {t, 2, ss}]; b4 =  Table[{t, 10}, {t, 2, ss + 2}];
Print[ff[x_] := 
 Which[{x[[1]], x[[2]]} == {-1, -1}, grap, x[[2]] == -1, x[[1]], 
 x[[1]] == -1, x[[2]], 1 == 1, Gr[{x[[1]], x[[2]], 0}]]; 
 Grid[Table[ff[{n, m}], {n, -1, ss}, {m, -1, ss}], Frame -> All, 
 Background -> {None, None, 
 Join[Table[b3[[s]] -> GrayLevel[gg], {s, 1, Length[b3]}],
 Table[b2[[s]] -> GrayLevel[gg], {s, 1, Length[b2]}]]}]];
Print[ff[x_] := 
 Which[{x[[1]], x[[2]]} == {-1, -1}, grap, x[[2]] == -1, x[[1]], 
 x[[1]] == -1, x[[2]], 1 == 1, Gr[{x[[1]], x[[2]], 1}]]; 
 Grid[Table[ff[{n, m}], {n, -1, ss}, {m, -1, ss}], Frame -> All, 
 Background -> {None, None, 
 Join[Table[b3[[s]] -> GrayLevel[gg], {s, 1, Length[b3]}],
 Table[b2[[s]] -> GrayLevel[gg], {s, 1, Length[b2]}]]}]]
	\end{verbatim}
\end{small}

\subsection{Calculation by CGSuite}
In this subsection, we give Tables of Grundy numbers using CGSuite.

\begin{figure}[H]
\begin{tabular}{cc}
\begin{minipage}[t]{0.48\textwidth}
\begin{center}
\includegraphics[height=2.7cm]{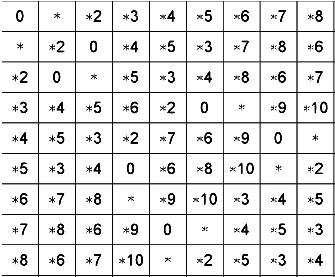}
\caption{Set $A$}
\label{nopass3}
\end{center}
\end{minipage}
\begin{minipage}[t]{0.48\textwidth}
\begin{center}
\includegraphics[height=2.7cm]{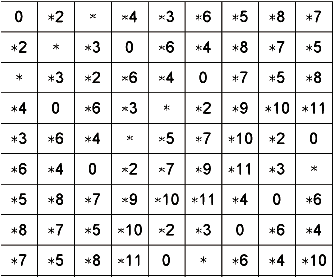}
\caption{Set $B$ }
\label{queenpass152}
\end{center}
\end{minipage}
\end{tabular}
\end{figure}

\subsection{Calculation by Python}
In this subsection, we give Tables of Grundy numbers using Python.
\begin{figure}[H]
\begin{tabular}{ccc}
\begin{minipage}[t]{0.33\textwidth}
\begin{center}
\includegraphics[height=2.9cm]{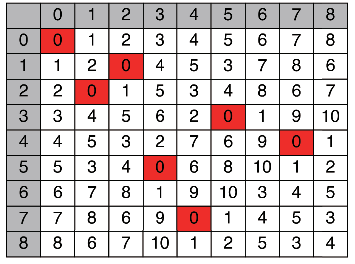}
\caption{Set $A$}
\label{nopass2}
\end{center}
\end{minipage}
\begin{minipage}[t]{0.33\textwidth}
\begin{center}
\includegraphics[height=2.9cm]{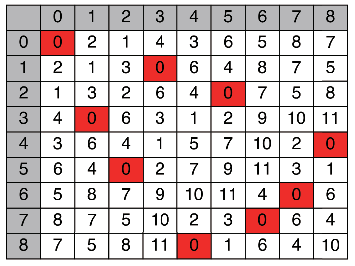}
\caption{Set $B$ }
\label{queenpass152b}
\end{center}
\end{minipage}
\begin{minipage}[t]{0.33\textwidth}
\begin{center}
\includegraphics[height=2.9cm]{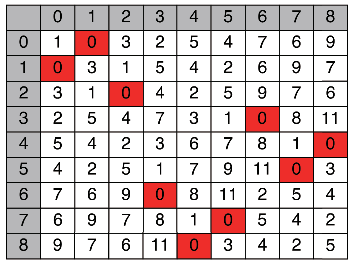}
\caption{Set $C$}
\label{queenplusone2}   
\end{center}
\end{minipage}
\end{tabular}
\end{figure}

\end{document}